# Gaussian-convolution-invariant shell approximation to spherically-symmetric functions

Alexandre G. Urzhumtsev[a,b] & Vladimir Y. Lunin[c]

[a]Université de Lorraine, Département de Physique, Faculté des Sciences et Technologies, BP 239, 54506 Vandoeuvre-les-Nancy, France

[b]Centre for Integrative Biology, Institut de Génétique et de Biologie Moléculaire et Cellulaire, CNRS–INSERM-UdS, 1 rue Laurent Fries, BP 10142, 67404 Illkirch, France

[c]Institute of Mathematical Problems of Biology RAS, Keldysh Institute of Applied Mathematics of Russian Academy of Sciences; Pushchino, 142290, Russia.

E-mail addresses: sacha@igbmc.fr (A.G. Urzhumtsev), lunin@impb.ru (V.Y. Lunin)

**Abstract**

We develop a class of functions $\Omega_N(\mathbf{x}; \mu, \nu)$ in $N$-dimensional space concentrated around a spherical shell of radius $\mu$ and such that, being convoluted with an isotropic Gaussian function, these functions do not change their expression but only a value of its 'width' parameter $\nu$. Isotropic Gaussian functions are a particular case of $\Omega_N(\mathbf{x}; \mu, \nu)$ corresponding to $\mu = 0$. Due to their features, these functions are an efficient tool to build approximations to smooth and continuous spherically-symmetric functions including oscillating ones. Atomic images in limited-resolution maps of electron density, electrostatic scattering potential and other scalar fields studies in physics, chemistry, biology, and other natural sciences are examples of such functions. We give simple analytic expressions of $\Omega_N(\mathbf{x}; \mu, \nu)$ for $N = 1, 2, 3$ and analyze properties of these functions. Representation of oscillating functions by a sum of $\Omega_N(\mathbf{x}; \mu, \nu)$ allows calculating distorted maps for the same cost as the respective theoretical fields. Using the chain rule and analytic expressions of the $\Omega_N(\mathbf{x}; \mu, \nu)$ derivatives makes simple refinement of parameters of the models which describe these fields.

**Keywords**

Spherical-shell approximation
Oscillating functions
Convolution invariance
Interference function
Three-dimensional space
Scalar-field modeling

## 1. Introduction

In various studies in natural sciences, experimentally estimated multi-dimensional scalar physical fields are modeled by a sum of a large number of spherically-symmetric contributions of elements of the object under study. Some examples are electrostatic potential in physics, electron density and scattering electrostatic potential in structural biology and sky brilliance in astronomy. Parameters of such models are determined and refined from comparison of maps of the theoretical fields with the experimental ones. A large number of contributions (up to a million of atoms in structural biology) and a need in repeating calculation of the maps in the process of refinement of the model parameters require development of efficient algorithms to calculate these maps. Due to experiment's features, experimental maps always miss some high-resolution details. Therefore, for a proper comparison, the maps of the theoretical fields should simulate this effect. While the contributions usually are peaky functions when calculating the theoretical fields themselves, they (rather, their images) become oscillating functions when calculating the limited-resolution maps of these fields. Additionally, experimental maps are blurred due to a presence of noise. This noise is usually modelled by the convolution of the map, or of the respective contributions, with a Gaussian function.

Let a function $f(\mathbf{x})$ be defined in a multi-dimensional space $\mathbf{R}^N$. We call this function isotropic or spherically-symmetric if it depends only on the length $|\mathbf{x}|$. We note this by $f(\mathbf{x}) = \bar{f}(|\mathbf{x}|)$ where its radial component $\bar{f}(x)$ is defined for $x \geq 0$. To approximate oscillating isotropic functions in multi-dimensional space $\mathbf{R}^N$, we introduce a class of 'shell functions' $\Omega_N(\mathbf{x}; \mu, \nu)$, $\mu \geq 0, \nu > 0$, also isotropic [1]. These functions are defined by the convolution of the singular uniform distribution concentrated at the surface of the sphere of radius $\mu$ in $\mathbf{R}^N$ with the $N$-dimensional isotropic Gaussian function of the 'width' $\nu$.

The approximation

$$f(\mathbf{x}) \approx \sum_{m=1}^{M} \kappa_m \Omega_N(\mathbf{x}; \mu_m, \nu_m) \tag{1}$$

may be considered as a generalization of an approximation to functions in one-dimensional space by a sum of positively and negatively weighted Gaussian functions, with the position of a peak referred to by $\mu_m$ (e.g., [2-5]).

By their definition, shell functions $\Omega_N(\mathbf{x}; \mu, \nu)$ are invariant with respect to the convolution with the Gaussian isotropic functions which just increases the parameter value $\nu$. As a consequence, approximation (1) is Gaussian-convolution-invariant in sense that its convolution with a Gaussian isotropic function changes only its inner parameters $\nu_m$. This feature is routinely used when $\mu_m = 0$, i.e., when each function $\Omega_N(\mathbf{x}; \mu_m, \nu_m)$ becomes a Gaussian one modeling a peaky contribution (e.g., [6-8]).

However, when modeling experimental maps of various experimentally obtained physical fields by sums of contributions, the contributions are often oscillating functions of the distance. This phenomenon is related to the Fourier series truncation and these oscillations have different names in different research fields. Such effect can be described by the convolution of the theoretical contribution with the interference function of the $N$-dimensional unit sphere. A decomposition of these standard interference functions into a sum of $\Omega_N(\mathbf{x}; \mu, \nu)$ results in an efficient tool to calculate limited-resolution maps of corresponding fields, even when this resolution varies from one map region to another, as it happens in cryo-sampled electron microscopy (cryoEM).

In this paper we present the mathematical background of the suggested approach. Some practical examples of the use of the developed technique in crystallography and in cryoEM may be found in [1,9,10].

Here and in what follows, to simplify equations, we omit $N = 1, 2, 3$ when the does not lead to a confusion, and use $x = |\mathbf{x}|$ for the length of the vector $\mathbf{x}$. Standard notations are used for the functions

$$\text{sinc}(x) = \frac{\sin(x)}{x}, \qquad \text{sinhc}(x) = \frac{\sinh(x)}{x}, \qquad \text{Si}(x) = \int_0^x \frac{\sin(t)}{t} dt \quad .$$

## 2. Modelling of map distortions

### 2.1. Normally distributed noise

For a proper numerical comparison of the maps of the theoretical fields with the experimental ones, the former should simulate distortions of the latter. One of the most common sources of distortions of the experimental maps like those of the electron density distribution is uncertainty in the atomic positions. This uncertainty can be dynamic if the atomic motion during the experiment cannot be neglected. Alternatively, these uncertainties may have a static origin, when averaging the measurements from a large number of samples, as in crystallography or in cryoEM. If such uncertainty in atomic coordinates is described by a probability distribution $P(\mathbf{x})$, a respective distorted map $f_{image}(\mathbf{x})$ may be modelled by the expected values. These values can be expressed by the convolution

$$f_{image}(\mathbf{x}) = \int_{\mathbf{R}^N} f(\mathbf{x} - \mathbf{u}) P(\mathbf{u}) dV_{\mathbf{u}} = f(\mathbf{x}) * P(\mathbf{x}) \quad . \tag{2}$$

Let both the contribution $f(\mathbf{x}; v) = a g_N(\mathbf{x}; v)$ to the theoretical map and the uncertainty $P(\mathbf{x}; v_0) = g_N(\mathbf{x}; v_0)$ be described by a normalized isotropic Gaussian function in $N$-dimensional space

$$g_N(\mathbf{x}; v) = \left(\frac{1}{2\pi v}\right)^{N/2} \exp\left(-\frac{|\mathbf{x}|^2}{2v}\right), \qquad \int_{\mathbf{R}^N} g_N(\mathbf{x}, v) dV_{\mathbf{x}} = 1 \quad . \tag{3}$$

Then the distorted image (2) of the contribution $f(\mathbf{x}; v)$ has the same form as the contribution itself but with a modified value of its parameter:

$$f_{image}(\mathbf{x}; v; v_0) = a g_N(\mathbf{x}; v + v_0) = f(\mathbf{x}; v + v_0) \quad . \tag{4}$$

This allows calculating the distorted map in the same way and for the same cost as the initial one. Inversely, values of the model parameters $a$ and $v$ can be easily estimated from the comparison of the calculated and experimental maps. The same property is generalized for the contributions $f(\mathbf{x})$ expressed by a sum of Gaussian isotropic functions.

### 2.2. Fourier transform of isotropic functions in multi-dimensional spaces

Another important source of the map distortion is loss of high-resolution data. Such distortion is described using the Fourier transform $\mathbf{F}(\mathbf{s}) = \mathcal{F}[f](\mathbf{s})$ of a function $f(\mathbf{x})$ in $N$-dimensional space $\mathbf{R}^N$. This transform is defined as

$$\mathbf{F}(\mathbf{s}) = \mathcal{F}[f](\mathbf{s}) = \int_{\mathbf{R}^N} f(\mathbf{x})\exp(i2\pi\mathbf{s}\cdot\mathbf{x})dV_\mathbf{x} \quad, \mathbf{s}\in\mathbf{R}^N \quad. \tag{6}$$

with the inverse Fourier transform

$$f(\mathbf{x}) = \mathcal{F}^{-1}[F](\mathbf{x}) = \int_{\mathbf{R}^N} \mathbf{F}(\mathbf{s})\exp(-i2\pi\mathbf{s}\cdot\mathbf{x})dV_\mathbf{s} \quad, \mathbf{x}\in\mathbf{R}^N \quad, \tag{5}$$

For a real-valued isotropic function $f(\mathbf{x})$, its Fourier transform $F(\mathbf{s})$ is also a real-valued isotropic function, and its direct and inverse Fourier transforms become the same. The radial components $\bar{f}(x)$ and $\bar{F}(s)$ of these functions are related to each other by a one-dimensional transform which expression depends on the dimension $N$:

$N = 1$

$$\bar{F}(s) = 2\int_0^\infty \bar{f}(x)\cos(2\pi s x)\, dx, s \geq 0 \quad ; \quad \bar{f}(x) = 2\int_0^\infty \bar{F}(s)\cos(2\pi s x)\, ds, x \geq 0 \quad, \tag{7}$$

$N = 2$

$$\bar{F}(s) = 2\pi\int_0^\infty x f(x) J_0(2\pi s x)\, dx, s \geq 0 \quad ; \quad \bar{f}(x) = 2\pi\int_0^\infty s F(s) J_0(2\pi s x)\, ds, x \geq 0 \quad, \tag{8}$$

$N = 3$

$$\bar{F}(s) = \frac{2}{s}\int_0^\infty x\bar{f}(x)\sin(2\pi s x)\, dx, s \geq 0 \quad, \bar{f}(x) = \frac{2}{x}\int_0^\infty s\bar{F}(s)\sin(2\pi s x)\, ds, x \geq 0 \quad, \tag{9}$$

Here and in what follows, $J_0(x)$ is the Bessel function of the zero order. Expressions (7)-(9) can be obtained by a straightforward integration using the parity of the function $f(x)$ for $N = 1$, integrating in polar coordinates for $N = 2$, and integrating in spherical coordinates for $N = 3$ with the polar axis chosen along the vector $\mathbf{s}$.

2.3.  *Loss of resolution*

Let a function $f(\mathbf{x})$ be defined in $N$-dimensional space $\mathbf{R}^N$, and $\mathbf{F}(\mathbf{s})$ be its Fourier transform. A map of $f(\mathbf{x})$ without high-resolution details can be formally described as the inverse Fourier transform

$$f_{image}(\mathbf{x}; d_0) = \int_{|\mathbf{s}|\leq d_0^{-1}} \mathbf{F}(\mathbf{s})\exp(-i2\pi\mathbf{s}\cdot\mathbf{x})dV_\mathbf{s}$$

of $\mathbf{F}(\mathbf{s})$ truncated to zero for $|\mathbf{s}| > d_0^{-1}$. Here $d_0$ is the resolution cut-off, or simply resolution of the map, and $\mathbf{s}\cdot\mathbf{x}$ is the scalar (dot) product of the vectors $\mathbf{s}$ and $\mathbf{x}$.

A straightforward numeric calculation of such map requires two Fourier transforms. Alternatively, using the convolution property of Fourier transform $\mathcal{F}^{-1}[\mathcal{F}[f]\cdot\mathcal{F}[g]] = f * g$, it can be defined, similar to (2), as a convolution

$$f_{image}(\mathbf{x}; d_0) = f(\mathbf{x}) * d_0^{-N} \pi_N(d_0^{-1}\mathbf{x}) \quad . \tag{10}$$

Here $\pi_N(\mathbf{x})$ is the inverse Fourier transform of the step-function in $N$-dimensional space

$$\Pi_N(\mathbf{s}) = \begin{cases} 1 & if \ |\mathbf{s}| \leq 1, \\ 0 & if \ |\mathbf{s}| > 1, \end{cases} \tag{11}$$

$$\pi_N(\mathbf{x}) = \mathcal{F}^{-1}[\Pi_M](\mathbf{x}) = \int_{\substack{\mathbf{s} \in \mathbf{R}^N \\ |\mathbf{s}| \leq 1}} \exp(-i2\pi \mathbf{s} \cdot \mathbf{x}) dV_\mathbf{s} \ , \mathbf{x} \in \mathbf{R}^N \ .$$

Functions $\pi_N(\mathbf{x})$, which are the interference functions of the unit sphere, or simply interference functions in what follows, are real-valued and isotropic. They satisfy conditions

$$\int_{\mathbf{R}^N} \pi_N(\mathbf{x}) dV_\mathbf{x} = 1, \qquad \pi_N(\mathbf{0}) = V_N = \int_{\substack{\mathbf{s} \in \mathbf{R}^N \\ |\mathbf{s}| \leq 1}} dV_\mathbf{s} \ ,$$

where $V_N$ is the volume of a unit sphere in $\mathbf{R}^N$.

For $N = 1, 2, 3$, the radial component of $\pi_N(\mathbf{x})$ can be expressed analytically as

$$\bar{\pi}_1(x) = 2\frac{\sin(2\pi x)}{2\pi x} \ ; \ \bar{\pi}_2(x) = 2\pi \frac{J_1(2\pi x)}{2\pi x} \ ; \ \bar{\pi}_3(x) = 4\pi \frac{\sin(2\pi x) - (2\pi x)\cos(2\pi x)}{(2\pi x)^3} \tag{12}$$

(**Fig. 1**). Here $J_1(x)$ is the Bessel function of the first order. Expressions for $\pi_N(\mathbf{x})$ are obtained applying (9) to (11) and using identity 6.561.5 from [11] for $N = 2$.

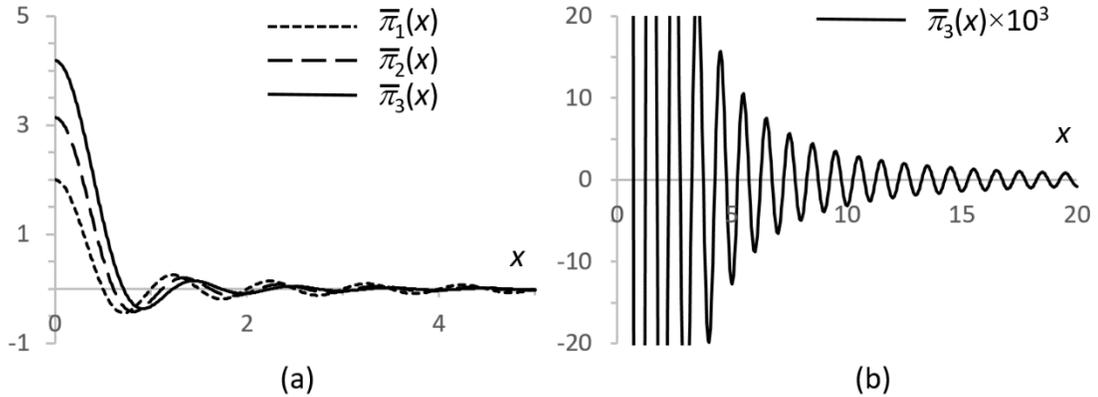

**Fig. 1.** a) Radial components of the functions $\pi_1(\mathbf{x}), \pi_2(\mathbf{x}), \pi_3(\mathbf{x})$; b) Zoom on the radial component of the function $\pi_3(\mathbf{x})$.

### 2.4. Combining the sources of image distortions

The two types of distortions defined above may be combined as

$$f_{image}(\mathbf{x}; v_0, d_0) = f(\mathbf{x}) * g_N(\mathbf{x}; v_0) * d_0^{-N} \pi_N(d_0^{-1}\mathbf{x})$$

where the convolution can be calculated in any order. Differently from (2) – (4), the result of convolution (10) does not have an analytic form when $f(\mathbf{x})$ is an isotropic Gaussian function. To overcome this obstacle, we approximate the interference functions $\pi_N(\mathbf{x})$ by the sums of specially designed shell functions $\Omega_N(\mathbf{x}; \mu, \nu)$ (see §3 and §4). These latter are defined as the convolution of the $N$-dimensional isotropic Gaussian function $g_N(\mathbf{x}; \nu)$ with the normalized uniform distribution $\delta_N(\mathbf{x}; \mu)$ on the surface of the sphere of the radius $\mu$ in $N$-dimensional space:

$$\Omega_N(\mathbf{x}; \mu, \nu) = \delta_N(\mathbf{x}; \mu) * g_N(\mathbf{x}; \nu) \quad .$$

The Gaussian function is a particular case of such functions,

$$\Omega_N(\mathbf{x}; 0, \nu) = g_N(\mathbf{x}; \nu) \quad .$$

By its definition and property (4), the class of such shell functions is invariant with respect to the convolution with the isotropic Gaussian functions

$$\Omega_N(\mathbf{x}; \mu, \nu) * g_N(\mathbf{x}; \nu_0) = \Omega_N(\mathbf{x}; \mu, \nu + \nu_0)$$

Below, we give analytical expressions for these functions $\Omega_N(\mathbf{x}; \mu, \nu)$ in one-, two- and three-dimensional spaces and analyze their principal features. We conclude our analysis by several practical applications of $\Omega_3(\mathbf{x}; \mu, \nu)$.

## 3. Shell-functions

### 3.1. Definition of the shell functions

We introduce a shell-function $\Omega_N(\mathbf{x}; \mu, \nu)$, $\mu \geq 0, \nu > 0$ in $\mathbf{R}^N$ as the convolution of a singular uniform distribution concentrated at the spherical surface of the radius $\mu$ and the normalized isotropic Gaussian function $g_N(\mathbf{x}; \nu)$, defined by (3). To obtain the explicit analytic expression for $\Omega_N(\mathbf{x}; \mu, \nu)$, we consider a family of characteristic functions of the spherical shells

$$p_N(\mathbf{x}; \mu, \varepsilon) = \frac{1}{S_N(\mu)\varepsilon} \begin{cases} 0 \text{ if } |\mathbf{x}| \leq \mu \\ 1 \text{ if } \mu < |\mathbf{x}| \leq \mu + \varepsilon \\ 0 \text{ if } |\mathbf{x}| > \mu + \varepsilon \end{cases}$$

approximating the singular distribution and normalized asymptotically

$$\lim_{\varepsilon \to 0} \int_{\mathbf{R}^N} p_N(\mathbf{r}; \mu, \varepsilon) dV_\mathbf{r} = 1 \quad .$$

Here $S_N(\mu)$ is the surface area of a sphere of the radius $\mu$ in $\mathbf{R}^N$. The Fourier transform $P_N(\mathbf{s}; \mu, \varepsilon) = \mathcal{F}[p_N]$ of the function $p_N(\mathbf{x}; \mu, \varepsilon)$ may be expressed as

$$P_N(\mathbf{s}; \mu, \varepsilon) = \frac{1}{S_N(\mu)\varepsilon} [H_N(\mathbf{s}; \mu + \varepsilon) - H_N(\mathbf{s}; \mu)] \tag{13}$$

where

$$H_N(\mathbf{s}; \mu) = \int_{\substack{\mathbf{x} \in \mathbf{R}^N \\ |\mathbf{x}| \leq \mu}} \exp(i 2\pi \mathbf{s} \cdot \mathbf{x}) dV_\mathbf{x} = \mu^N \pi_N(\mu \mathbf{s}) \quad ,$$

and $\pi_N(\mathbf{x})$ is the interference function defined by (11). The limit of $P_N(\mathbf{s}; \mu, \varepsilon)$ is

$$P_N^0(\mathbf{s}; \mu) = \lim_{\varepsilon \to 0} P_N(\mathbf{s}; \mu, \varepsilon) = \frac{1}{S_N(\mu)} \lim_{\varepsilon \to 0} \frac{H_N(\mathbf{s}; \mu + \varepsilon) - H_N(\mathbf{s}; \mu)}{\varepsilon} = \frac{1}{S_N(\mu)} \frac{\partial (\mu^N \pi_N(\mu \mathbf{s}))}{\partial \mu}$$

With (12), the radial part $\overline{P}_N^0(s; \mu)$ for $N = 1, 2, 3$ can be expressed analytically as

$$\bar{P}_1^0(s;\mu) = \cos(2\pi\mu s) \ ; \quad \bar{P}_2^0(s;\mu) = J_0(2\pi\mu s) \ ; \quad \bar{P}_3^0(s;\mu) = \frac{\sin(2\pi\mu s)}{2\pi\mu s}$$

(with 8.472.3 from [11] for $N = 2$). Using the convolution theorem for the Fourier transform, we calculate $\Omega_N(\mathbf{x};\mu,\nu)$ as

$$\Omega_N(\mathbf{x};\mu,\nu) = \mathcal{F}^{-1}[\bar{P}_N^0(\mathbf{s};\mu) \cdot G_N(\mathbf{s};\nu)] \ , \tag{14}$$

where

$$G_N(\mathbf{s};\nu) = \exp(-2\pi^2\nu|\mathbf{s}|^2) \tag{15}$$

is the Fourier transform of the Gaussian function (3).

For the spherically symmetric functions $\bar{P}_N^0(\mathbf{s};\mu)$ and $G_N(\mathbf{s};\nu)$, with $N = 1, 2, 3$, their inverse Fourier transform (14) can be calculated analytically giving the radial components equal to

$$\bar{\Omega}_1(x;\mu,\nu) = 2\int_0^\infty \cos(2\pi\mu s)\exp(-2\pi^2\nu s^2)\cos(2\pi s x)\,ds$$

$$= \frac{1}{2} \cdot \frac{1}{\sqrt{2\pi\nu}}\left[\exp\left(-\frac{(x-\mu)^2}{2\nu}\right) + \exp\left(-\frac{(x+\mu)^2}{2\nu}\right)\right] \tag{16}$$

$$\bar{\Omega}_2(x;\mu;\nu) = 2\pi\int_0^\infty J_0(2\pi\mu s)\exp(-2\pi^2\nu s^2)\,s\,J_0(2\pi s x)\,ds$$

$$= \frac{1}{2\pi\nu}\exp\left(-\frac{x^2+\mu^2}{2\nu}\right)I_0\left(\frac{x\mu}{\nu}\right) \tag{17}$$

$$\bar{\Omega}_3(x;\mu,\nu) = \frac{2}{x}\int_0^\infty \frac{\sin(2\pi\mu s)}{2\pi\mu s}\exp(-2\pi^2\nu s^2)\,s\sin(2\pi s x)\,ds$$

$$= \frac{1}{4\pi x\mu} \cdot \frac{1}{\sqrt{2\pi\nu}}\left[\exp\left(-\frac{(x-\mu)^2}{2\nu}\right) - \exp\left(-\frac{(x+\mu)^2}{2\nu}\right)\right] \tag{18}$$

(see 3.898.1, 3.898.2, 6.633.2 in [11]). Here $I_0(x)$ is the modified Bessel function of the zero order. **Fig. 2** shows $\kappa\bar{\Omega}_N(x;\mu,\nu)$ for several values of its parameters. Formulas (16) and (18) can be rewritten in the same form as (17)

$$\bar{\Omega}_1(x;\mu,\nu) = \frac{1}{(2\pi\nu)^{1/2}}\exp\left(-\frac{x^2+\mu^2}{2\nu}\right)\cosh\left[\frac{x\mu}{\nu}\right] \tag{19}$$

$$\bar{\Omega}_3(x;\mu,\nu) = \frac{1}{(2\pi\nu)^{3/2}}\exp\left(-\frac{x^2+\mu^2}{2\nu}\right)\sinh c\left(\frac{x\mu}{\nu}\right) \ , \tag{20}$$

more comprehensive but less appropriated for straightforward calculations with large values of $x\mu\nu^{-1}$. More practicable may be using expressions where the functions in the brackets have a limit equal to 1 when $x \to \infty$ :

$$\bar{\Omega}_1(x;\mu,\nu) = \frac{1}{\sqrt{2\pi\nu}}\exp\left(-\frac{(x-\mu)^2}{2\nu}\right)\cdot\frac{1}{2}\left[1 + \exp\left(-\frac{2x\mu}{\nu}\right)\right] \ ,$$

$$\bar{\Omega}_2(x;\mu;\nu) = \frac{1}{\sqrt{2\pi\nu}} \exp\left(-\frac{(x-\mu)^2}{2\nu}\right) \cdot \frac{1}{2\pi\sqrt{x\mu}} \left[\left(\frac{2\pi x\mu}{\nu}\right)^{\frac{1}{2}} \cdot I_0\left(\frac{x\mu}{\nu}\right) \exp\left(-\frac{x\mu}{\nu}\right)\right],$$

$$\bar{\Omega}_3(x;\mu,\nu) = \frac{1}{\sqrt{2\pi\nu}} \exp\left(-\frac{(x-\mu)^2}{2\nu}\right) \cdot \frac{1}{4\pi x\mu} \left[1 - \exp\left(-\frac{2x\mu}{\nu}\right)\right].$$

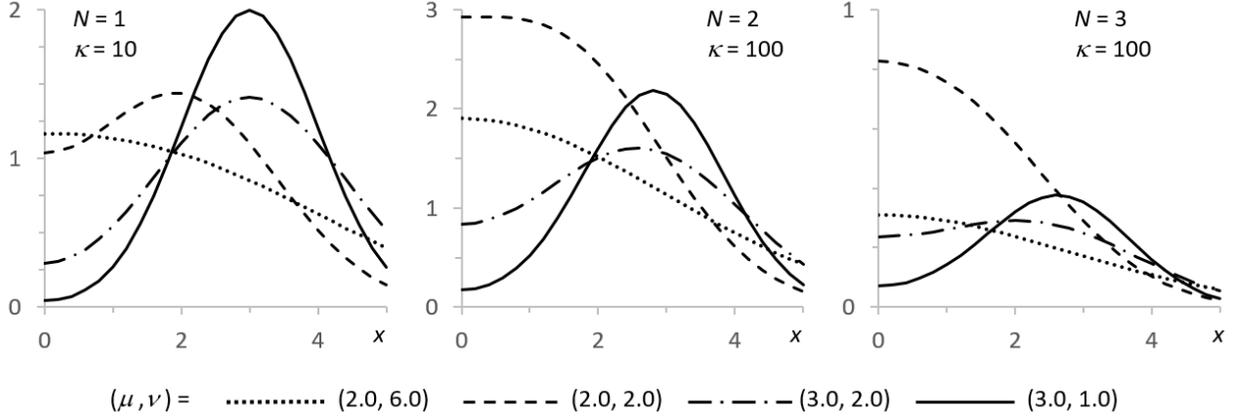

**Fig. 2.** Radial components of the functions $\kappa\bar{\Omega}_N(x;\mu,\nu)$ for various values of the parameters.

### 3.2. Properties of the shell functions

**Property 1.**

$$\lim_{x\to\infty} \bar{\Omega}_N(x;\mu,\nu) = 0 \; ; \; \bar{\Omega}_N(0;\mu,\nu) = \frac{1}{(2\pi\nu)^{N/2}} \exp\left(-\frac{\mu^2}{2\nu}\right)$$

**Property 2.**

$$\bar{\Omega}_N(x;\mu,\nu) = \frac{1}{(2\pi\nu)^{N/2}} \exp\left(-\frac{x^2+\mu^2}{2\nu}\right)\left[1 + O\left(\left(\frac{x\mu}{\nu}\right)^2\right)\right] \; for \; \frac{x\mu}{\nu} \ll 1$$

**Property 3.** Gaussian functions $g_N(\mathbf{x};\nu_0)$ are a particular case of the shell functions

$$g_N(\mathbf{x};\nu) = \Omega_N(\mathbf{x};0,\nu)$$

**Proofs.** Properties 1-3 follow directly from expressions (17), (19), (20).

**Property 4. Parameters rescaling.**

$$\bar{\Omega}_N\left(\frac{x}{\alpha};\mu,\nu\right) = \alpha^N \bar{\Omega}_N(x;\alpha\mu,\alpha^2\nu) \tag{21}$$

**Property 5. Gaussian-convolution independence.**

$$\Omega_N(\mathbf{x};\mu,\nu) * g_N(\mathbf{x};\nu_0) = \Omega_N(\mathbf{x};\mu,\nu+\nu_0)$$

**Proofs.** Properties 4 and 5 follow directly from the function definition (14) and the convolution

theorem (for Property 5):

**Property 6. Unimodality of the function.** For $\mu^2 \leq N\nu$, function $\bar{\Omega}_N(x;\mu,\nu)$ is a decreasing function on the demi-axis $x > 0$ with a maximum in $x = 0$. For $\mu^2 > N\nu$, function $\bar{\Omega}_M(x;\mu,\nu)$ has a local minimum in $x = 0$ and a single local extremum (maximum) on the demi-axis $x > 0$.

**Proof.** Essentially, this property means that function $\bar{\Omega}_N(x;\mu,\nu)$ cannot have simultaneously a peak on the demi-axis and at the origin. We propose the proof for practically interesting spaces with $N = 1, 2, 3$. After rescaling the variable and the parameter as

$$t = \frac{x\mu}{\nu} > 0 \,; a = \frac{\mu^2}{\nu} > 0$$

and ignoring all factors constant, functions $\bar{\Omega}_N(x;\mu,\nu)$ in (17), (19), (20) are reduced to

$$\omega_1(t;a) = \exp\left(-\frac{t^2}{2a}\right)\cosh(t)$$

$$\omega_2(t;a) = \exp\left(-\frac{t^2}{2a}\right)I_0(t)$$

$$\omega_3(t;a) = \frac{1}{t}\exp\left(-\frac{t^2}{2a}\right)\sinh(t)$$

These functions have an extremum when:

$$\frac{\partial \omega_1(t;a)}{\partial t} = \left[t\exp\left(-\frac{t^2}{2a}\right)\cosh(t)\right]\left[\frac{\tanh(t)}{t} - \frac{1}{a}\right] = 0$$

$$\frac{\partial \omega_2(t;a)}{\partial t} = \left[t I_0(t)\exp\left(-\frac{t^2}{2a}\right)\right]\left[\frac{I_1(t)}{t I_0(t)} - \frac{1}{a}\right] = 0$$

$$\frac{\partial \omega_3(t;a)}{\partial t} = \left[\sinh(t)\exp\left(-\frac{t^2}{2a}\right)\right]\left[\frac{t\cosh(t) - \sinh(t)}{t^2 \sinh(t)} - \frac{1}{a}\right] = 0$$

The functions inside the left brackets are always positive except at $t = 0$ where their value is equal to zero. The functions inside the right brackets

$$\chi_1(t) = \frac{\tanh(t)}{t} \,; \; \chi_2(t) = \frac{I_1(t)}{t I_0(t)} \,; \; \chi_3(t) = \frac{t\cosh(t) - \sinh(t)}{t^2 \sinh(t)}$$

are monotonously decreasing with the asymptote $\chi_N(t) = t^{-1}$ (**Fig. 3**). Their values in the origin are $\chi_N(0) = N^{-1}$ which, with $a = \mu^2\nu^{-1}$, proves the property.

**Property 7.**
Function $\bar{\Omega}_N(x;\mu,\nu)$ is decreasing in $x = \mu$.
**Proof.**
We propose the proof for practically interesting spaces $N = 1, 2, 3$ using known inequalities

$$\tanh(a) < a \,; \; I_1(a) < I_0(a) \,; \; \tanh(a) > \frac{a}{a+1}$$

Applying them to the rescaled derivatives $\omega_N(t;a)$ of the functions $\bar{\Omega}_N(x;\mu,\nu)$ we get

$$\frac{\partial \omega_N(a;a)}{\partial t} < 0, N = 1, 2, 3$$

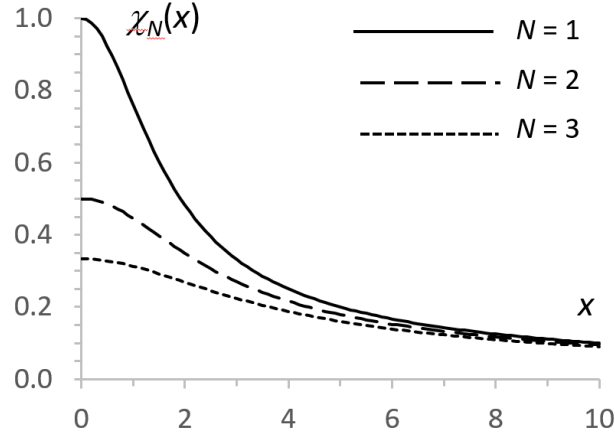

**Fig. 3.** Functions $\chi_N(x)$.

## 4. Decomposition of the interference functions

Algorithms and programs to build the shell-approximation to the functions in three-dimensional space are described in [12]. To build an approximation (1) to a given function, one chooses the value of the decomposition interval, $|\mathbf{x}| \leq x_{max}$, and the required accuracy. The parameter values $(\mu_m, \nu_m, \kappa_m)$ are estimated from the position and shape of the respective ripples of the function at this interval. When a higher accuracy is required, extra terms are added modeling the residual peaks. The peak in the origin may be modeled by several terms including ones with $\mu_m > 0$. These initial estimates are then optimized using the standard minimizer L-BFGS [13].

**Table 1** shows the coefficients $(\mu_m, \nu_m, \kappa_m)$ of the approximation to the interference function

$$\pi_3(\mathbf{x}) \approx \sum_{m=1}^{M} \kappa_m \Omega_3(\mathbf{x}; \mu_m, \nu_m)$$

at the interval $|\mathbf{x}| \leq 20$, one term per ripple, positive or negative, with 40 terms in total, giving the maximal point-by-point discrepancy $6.0 \times 10^{-4}$.

**Table 1.** Coefficients of the shell approximation to the unit-sphere interference function $\pi_3(\mathbf{x})$ in three-dimensional space calculated for $|\mathbf{x}| \leq 20$.

| m | $\mu$ | $\nu$ | $\kappa$ | m | $\mu$ | $\nu$ | $\kappa$ |
|---|---|---|---|---|---|---|---|
| 1 | 0.000000 | 0.149577 | 3.974732 | 21 | 10.493471 | 0.016800 | 1.338413 |
| 2 | 0.839261 | 0.081904 | -4.722464 | 22 | 10.993226 | 0.016872 | -1.336079 |
| 3 | 1.411624 | 0.054497 | 2.987349 | 23 | 11.493476 | 0.016796 | 1.334134 |
| 4 | 1.952697 | 0.042378 | -2.323489 | 24 | 11.994229 | 0.016080 | -1.332702 |
| 5 | 2.464160 | 0.036608 | 2.056706 | 25 | 12.495735 | 0.018301 | 1.332798 |
| 6 | 2.974455 | 0.031845 | -1.882625 | 26 | 12.996825 | 0.019754 | -1.331953 |
| 7 | 3.483748 | 0.030437 | 1.823866 | 27 | 13.497168 | 0.018059 | 1.330944 |
| 8 | 3.980074 | 0.028966 | -1.763604 | 28 | 13.997713 | 0.016041 | -1.329831 |
| 9 | 4.483824 | 0.025723 | 1.660668 | 29 | 14.498717 | 0.016684 | 1.328876 |
| 10 | 4.987578 | 0.024247 | -1.611679 | 30 | 14.999067 | 0.016718 | -1.328008 |
| 11 | 5.485233 | 0.021983 | 1.532103 | 31 | 15.496648 | 0.016876 | 1.337655 |
| 12 | 5.992215 | 0.020754 | -1.482867 | 32 | 15.996165 | 0.016851 | -1.337766 |
| 13 | 6.490989 | 0.020406 | 1.490286 | 33 | 16.496408 | 0.016843 | 1.337641 |
| 14 | 6.991365 | 0.019200 | -1.425480 | 34 | 16.996653 | 0.016761 | -1.337527 |
| 15 | 7.492524 | 0.018104 | 1.392489 | 35 | 17.496881 | 0.016769 | 1.337448 |
| 16 | 7.991915 | 0.017389 | -1.356501 | 36 | 17.997110 | 0.016881 | -1.337390 |
| 17 | 8.493796 | 0.017035 | 1.347738 | 37 | 18.497337 | 0.017033 | 1.337344 |
| 18 | 8.994394 | 0.016775 | -1.347498 | 38 | 18.997558 | 0.017155 | -1.337322 |
| 19 | 9.494283 | 0.016393 | 1.325074 | 39 | 19.497772 | 0.017199 | 1.337314 |
| 20 | 9.994853 | 0.016184 | -1.304771 | 40 | 19.997212 | 0.016972 | -1.337255 |

A similar procedure has been used to calculate approximations

$$\pi_N(\mathbf{x}) \approx \sum_{m=1}^{M} \kappa_m \Omega_N(\mathbf{x}; \mu_m, \nu_m)$$

for $N = 1, 2$, substituting $\Omega_3(\mathbf{x}; \mu, \nu)$ by $\Omega_1(\mathbf{x}; \mu, \nu)$ and by $\Omega_2(\mathbf{x}; \mu, \nu)$, respectively. **Table 2** and **Table 3** give the coefficients of such approximations, with in the maximal error equal to $7.8 \times 10^{-5}$ and to $4.1 \times 10^{-4}$, respectively. **Fig. 4** illustrates the discrepancy, as a function of $|\mathbf{x}|$, at the interval $|\mathbf{x}| \leq 10$ between the $\pi_N(\mathbf{x})$ and shell-approximations to them with 21 terms given in **Tables 1-3**.

**Table 2.** Coefficients of the shell approximation to the one-dimensional unit-sphere interference function $\pi_1(\mathbf{x})$ calculated for $|\mathbf{x}| \leq 10$.

| m | $\mu$ | $\nu$ | $\kappa$ | m | $\mu$ | $\nu$ | $\kappa$ |
|---|---|---|---|---|---|---|---|
| 1 | 0.000000 | 0.131571 | 2.173815 | 11 | 5.240482 | 0.047119 | 0.076789 |
| 2 | 0.591792 | 0.107565 | -1.638004 | 12 | 5.740105 | 0.046499 | -0.069314 |
| 3 | 1.148469 | 0.082861 | 0.680768 | 13 | 6.241548 | 0.045711 | 0.062843 |
| 4 | 1.687606 | 0.065512 | -0.343475 | 14 | 6.742026 | 0.045170 | -0.057479 |
| 5 | 2.211743 | 0.055550 | 0.215743 | 15 | 7.241991 | 0.044328 | 0.052689 |
| 6 | 2.725149 | 0.050571 | -0.158344 | 16 | 7.742989 | 0.043498 | -0.048486 |
| 7 | 3.232188 | 0.048228 | 0.127375 | 17 | 8.244333 | 0.042899 | 0.045129 |
| 8 | 3.736620 | 0.047272 | -0.108157 | 18 | 8.747921 | 0.043574 | -0.043019 |
| 9 | 4.239295 | 0.047227 | 0.095144 | 19 | 9.242462 | 0.044127 | 0.040683 |
| 10 | 4.739889 | 0.047297 | -0.085197 | 20 | 9.732445 | 0.039877 | -0.035737 |
|   |   |   |   | 21 | 10.230986 | 0.031269 | 0.030208 |

**Table 3.** Coefficients of the shell approximation to the two-dimensional unit-sphere interference function $\pi_2(\mathbf{x})$ calculated for $|\mathbf{x}| \leq 10$.

| m | $\mu$ | $\nu$ | $\kappa$ | m | $\mu$ | $\nu$ | $\kappa$ |
|---|---|---|---|---|---|---|---|
| 1 | 0.000000 | 0.121565 | 2.488048 | 11 | 5.366582 | 0.032426 | 0.408413 |
| 2 | 0.742617 | 0.073698 | -2.253003 | 12 | 5.866127 | 0.031651 | -0.385374 |
| 3 | 1.304418 | 0.054291 | 1.267576 | 13 | 6.367222 | 0.030875 | 0.364248 |
| 4 | 1.833737 | 0.044795 | -0.892266 | 14 | 6.867885 | 0.030196 | -0.345915 |
| 5 | 2.347657 | 0.040104 | 0.718677 | 15 | 7.368364 | 0.029513 | 0.329532 |
| 6 | 2.854680 | 0.037505 | -0.619125 | 16 | 7.868835 | 0.028856 | -0.314738 |
| 7 | 3.358746 | 0.035779 | 0.552030 | 17 | 8.369306 | 0.028241 | 0.301427 |
| 8 | 3.861382 | 0.034511 | -0.502009 | 18 | 8.869483 | 0.027629 | -0.289143 |
| 9 | 4.362868 | 0.033351 | 0.462110 | 19 | 9.368906 | 0.026759 | 0.276422 |
| 10 | 4.865655 | 0.032578 | -0.431078 | 20 | 9.867940 | 0.025535 | -0.261959 |
|  |  |  |  | 21 | 10.366720 | 0.023456 | 0.205540 |

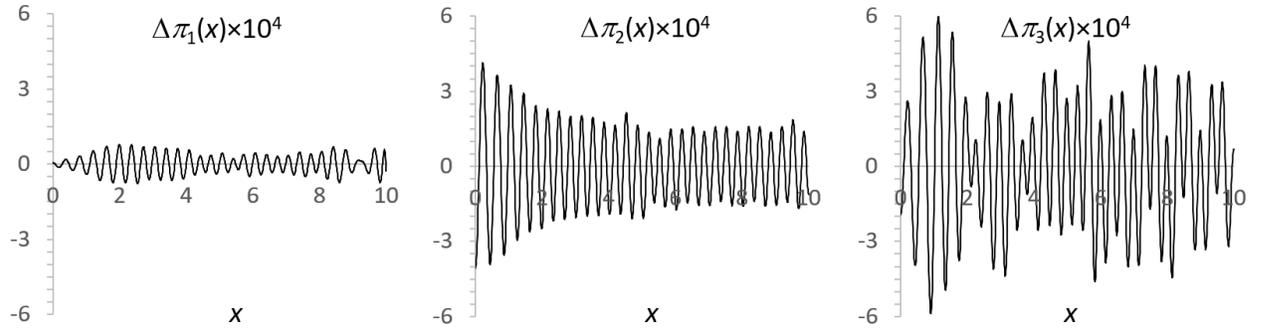

**Fig. 4.** Discrepancy between the radial components of the interference functions $\pi_N(\mathbf{x})$ and the spherical-shell approximations to them by a sum of 21 shell functions for $|\mathbf{x}| \leq 10$.

## 5. Examples of the shell decomposition in three-dimensional space

### 5.1. Images of the Gaussian and Gaussian-like function

Diffraction methods such as X-ray crystallography, cryoEM, microcrystal electron diffraction, X-FEL and others determine atomic structure studying spatial scalar functions, *e.g.*, an electron density distribution, which are sums of atomic contributions. In structural biology, an isotropic Gaussian atomic disorder is expressed by so-called atomic displacement parameter defined as $B = 8\pi^2 \nu$. Atomic contributions, such as an atomic electron density distribution and an electrostatic scattering potential, are traditionally approximated by a linear combination of a few three-dimensional Gaussians [6-8]

$$f(\mathbf{x}) \approx \sum_{k=1}^{K^{Gauss}} a^{(k)} g_3\left(\mathbf{x}; \frac{B^{(k)}}{8\pi^2}\right), \qquad \mathbf{x} \in \mathbf{R}^3 \qquad . \qquad (22)$$

(Actually, the Gaussian approximations are built rather to the scattering functions $\bar{F}(s)$, which are the Fourier transforms of the contribution of immobile atoms.)

According to (10) and (21), the image of the Gaussian function $f(\mathbf{x}) = a g_3(\mathbf{x}; b)$ at the

resolution $d_0$ is

$$f^{(d)}(\mathbf{x}; d_0) = ag_3(\mathbf{x}; b) * d_0^{-N} \pi_N(d_0^{-1}\mathbf{x}) \approx ag_3(\mathbf{x}; b) * d_0^{-3} \sum_{m=1}^{M} \kappa_m \Omega_3(d_0^{-1}\mathbf{x}; \mu_m, \nu_m)$$

$$= a \sum_{m=1}^{M} \kappa_m \Omega_3(\mathbf{x}; d_0\mu_m, b + d_0^2 \nu_m).$$

If an additional harmonic noise characterized by the parameter $\nu_0$ is introduced, the image above is updated simply as

$$f_{image}(\mathbf{x}; \nu_0, d_0) = a \sum_{m=1}^{M} \kappa_m \Omega_3(\mathbf{x}; d_0\mu_m, b + d_0^2 \nu_m + \nu_0) \quad.$$

Then, an image of any atomic contribution (22) at any position $\mathbf{x}_n$ with any displacement parameter $B_n$ and at any resolution $d_0$ can be expressed analytically as

$$f^{(d)}(\mathbf{x}; B_n, d_0) \approx \sum_{k=1}^{K^{Gauss}} a^{(k)} \sum_{m=1}^{M} \kappa_m \Omega_3 \left( \mathbf{x} - \mathbf{x}_n; \mu_m d_0, \nu_m d_0^2 + \frac{B_n + B^{(k)}}{8\pi^2} \right). \qquad (23)$$

Alternatively, when the resolution value $d_0$ is known and fixed, one can simplify (23). First, an image of an immobile atom of the given chemical type at the resolution $d_0$ is calculated as the inverse Fourier transform (9) of its scattering factor $\bar{F}(s)$ cut out for $|\mathbf{s}| > d_0^{-1}$. Then one can build a respective approximation

$$f^{(d)}(\mathbf{x}; 0, d_0) \approx \sum_{m=1}^{M'} K^{(m)} \Omega_3 \left( \mathbf{x}; R^{(m)}, \frac{B^{(m)}}{8\pi^2} \right) \qquad (24)$$

directly to this image resulting in

$$f^{(d)}(\mathbf{x}; B_n, d_0) \approx \sum_{m=1}^{M'} K^{(m)} \Omega_3 \left( \mathbf{x}; R^{(m)}, \frac{B_n + B^{(m)}}{8\pi^2} \right).$$

This reduces roughly by $K^{Gauss}$ times the number of terms of the approximation in comparison with (23) and therefore the CPU time required to calculate, but requires to calculate preliminary the series coefficients (24) for all types of atoms present in the structure [10].

*5.2. Images of non-Gaussian-like functions*

The Coulomb or gravitational potential are defined as $u_0(\mathbf{x}) = K|\mathbf{x}|^{-1}$ where coefficient $K$ includes both the factor with the universal physical constants and the charge (or mass) of the source of the potential. These functions cannot be presented in the multi-Gaussian form (22). However, one can express in this form their images introducing artificially some small value $\nu$ of the uncertainty (see Appendix A). Then the resulted function

$$u(\mathbf{x}; \nu) = (K|\mathbf{x}|^{-1}) \operatorname{erf}\left[|\mathbf{x}|(2\nu)^{-1/2}\right]$$

becomes representative by a sum (22) as **Fig. 5a** shows.

Alternatively, if the calculations supposed to be carried out at a known resolution $d_0$, one can calculate first the image of this potential for an 'immobile charge' at this resolution using the known Fourier transform $U_0(\mathbf{s}) = K/(\pi|\mathbf{s}|^2)$ of this potential (see Appendix A). According to (9), this image of this spherically symmetric function is

$$\bar{u}^{(d)}(x; d_0) = \frac{2}{x} \int_0^{d_0^{-1}} s \frac{K}{\pi s^2} \sin(2\pi x s) \, ds = 4K d_0^{-1} \frac{\text{Si}(2\pi x d_0^{-1})}{2\pi x d_0^{-1}} \quad . \tag{25}$$

Function (25) is not oscillating (**Fig. 5b**) but can be also accurately approximated by a series (1) of $\Omega_3(\mathbf{x}; \mu_m, \nu_m)$ terms (**Table 4**). **Figs. 5b-5d** illustrate the results of such approximation to the function $|\mathbf{x}|^{-1} \text{Si}(|\mathbf{x}|)$, $N = 3$, in the interval $|\mathbf{x}| \leq 10$, by a series with $M = 4$ first terms (maximal approximation error at his interval is equal to $1.5 \times 10^{-2}$), with $M = 16$ first terms (maximal error is equal to $7.7 \times 10^{-4}$) and with $M = 33$ terms (maximal error $2.2 \times 10^{-4}$).

**Table 4.** Coefficients of the shell approximation to the function $\text{Si}(|\mathbf{x}|)/|\mathbf{x}|$ in three-dimensional space calculated for $|\mathbf{x}| \leq 10$ (**Figs. 5b-d**). Note importance of the term with $m = 4$ compensating the 'tail' of the term $m = 1$.

| m | μ | ν | κ | m | μ | ν | κ |
|---|---|---|---|---|---|---|---|
| 1 | 0.000000 | 6.197172 | 268.311753 | 17 | 0.140120 | 0.003620 | -0.000080 |
| 2 | 0.000000 | 0.163080 | 0.806227 | 18 | 0.370197 | 0.002062 | -0.000136 |
| 3 | 2.444275 | 3.331012 | -206.872124 | 19 | 0.559895 | 0.004527 | 0.001692 |
| 4 | 10.917662 | 14.020006 | 352.636534 | 20 | 0.839973 | 0.002748 | 0.001209 |
| 5 | 0.020011 | 0.006256 | -0.000039 | 21 | 1.040030 | 0.005487 | -0.007025 |
| 6 | 0.420455 | 0.011020 | 0.005567 | 22 | 1.350049 | 0.003893 | -0.003426 |
| 7 | 0.871281 | 0.013373 | -0.032196 | 23 | 1.580375 | 0.005814 | 0.016004 |
| 8 | 1.391771 | 0.021336 | 0.138810 | 24 | 1.960182 | 0.004838 | 0.015088 |
| 9 | 1.979248 | 0.021204 | -0.245133 | 25 | 2.210897 | 0.006995 | -0.038062 |
| 10 | 2.560773 | 0.010795 | 0.050606 | 26 | 2.990484 | 0.013585 | -0.041777 |
| 11 | 3.482397 | 0.048566 | 0.448946 | 27 | 4.460091 | 0.025375 | 0.099499 |
| 12 | 5.030498 | 0.057977 | -0.636492 | 28 | 5.439990 | 0.037814 | -0.099893 |
| 13 | 5.900254 | 0.054990 | -0.566002 | 29 | 6.610000 | 0.050331 | 0.249889 |
| 14 | 7.510280 | 0.104842 | 1.231657 | 30 | 9.109993 | 0.037952 | -0.299975 |
| 15 | 8.370087 | 0.054051 | 0.632798 | 31 | 2.729991 | 0.002894 | 0.002875 |
| 16 | 10.100056 | 0.182008 | -3.032464 | 32 | 3.709994 | 0.002962 | 0.004538 |
|   |   |   |   | 33 | 4.049998 | 0.002976 | -0.002545 |

## 6. Discussion

The designed functions $\Omega_N(\mathbf{x}; \mu; \nu)$ in multi-dimensional spaces have a number of useful features making them an efficient tool in solution of numerous practical problems, in particular in three-dimensional space. Invariance of the class of these functions with respect to a convolution with the Gaussian functions is one of these features. Approximation to the interference functions of the unit sphere by a series of $\Omega_N(\mathbf{x}; \mu; \nu)$ is an example of such tools for modeling variable-resolution maps of various scalar fields in two- and three-dimensional spaces. Coefficients of these approximation to the interference functions are given in this work allowing a direct use of such modeling. The programs to calculate approximations by $\Omega_N(\mathbf{x}; \mu; \nu)$ are available by request from the authors (AU).

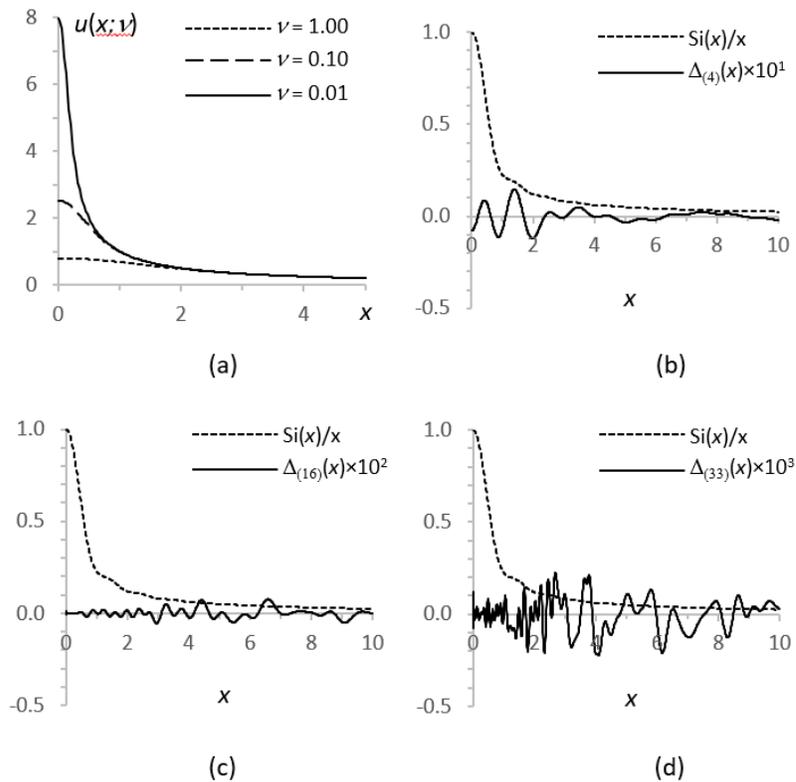

**Figure 5.** a) Radial component of the Coulomb potential with the introduced harmonic uncertainty $\nu$, the value indicated in the legend; b-d) Radial component of the function $|\mathbf{x}|^{-1}\,\text{Si}(|\mathbf{x}|)$, and the discrepancy $\Delta(x; M)$ between this function and the spherical-shell approximation to it calculated with $M = 4$, 16 and 33 terms.

**Acknowledgment**

AU acknowledges the French Infrastructure for Integrated Structural Biology FRISBI [ANR-10-INBS-05] and the support and the use of resources of Instruct-ERIC through the R&D pilot scheme APPID 2683.

**Appendix A. Coulomb potential and its images**

A well-known procedure to calculate the Fourier transform of the Coulomb (or gravitational) potential in three-dimensional space $u_0(\mathbf{x}) = K|\mathbf{x}|^{-1}$ consists in calculating it initially for a more general function, known as the Yukawa potential, $u(\mathbf{x}; \lambda) = K \exp(-\lambda|\mathbf{x}|) |\mathbf{x}|^{-1}$, with $\lambda > 0$.

Using the spherical symmetry of the function, we calculate the transform as (9) along the axis $\mathbf{s} = s\mathbf{O}_z$ which gives

$$\bar{U}_\lambda(s) = \frac{4\pi K}{\lambda^2 + 4\pi^2 s^2}$$

Its limit for $\lambda \to 0$ defines the Fourier transform of the Coulomb potential as $\bar{U}_0(s) = K(\pi s^2)^{-1}$.

The image of this potential in presence of a normally-distributed uncertainty in the source position is described by the convolution giving ([11], 3.952.6[3])

$$u(\mathbf{x}; \nu) = \frac{K}{|\mathbf{x}|} * g_3(\mathbf{x}; \nu) = \mathcal{F}^{-1}\left[\frac{K}{\pi|\mathbf{s}|^2} \cdot \exp(-2\pi^2 \nu |\mathbf{s}|^2)\right] = \frac{K}{|\mathbf{x}|} \operatorname{erf}\left[\frac{|\mathbf{x}|}{(2\nu)^{1/2}}\right]$$

Here the right-hand expression can be obtained using (9) and identity 3.952.6[3] in [11].